\documentclass[a4paper,11pt,reqno]{amsart}

\usepackage{amsfonts}
\usepackage{amssymb}
\usepackage{amscd}
\usepackage{palatino}
\usepackage{amsthm}
\usepackage{a4wide}
\usepackage{mathrsfs}

\usepackage[centertags]{amsmath}
\usepackage{eucal,dsfont,accents,bbm,amsfonts,amssymb,amsthm,amsopn}

\usepackage[usenames]{color}
\definecolor{citegreen}{rgb}{0,0.6,0}
\definecolor{refred}{rgb}{0.8,0,0}
\usepackage[colorlinks, citecolor=citegreen, linkcolor=refred]{hyperref}

\usepackage[]{draftcopy}
\draftcopyName{ Draft}{240}
\draftcopySetScale{270}

\theoremstyle{plain}

\newtheorem{teo}{Theorem}[section]
\newtheorem{lemma}[teo]{Lemma}

\newtheorem{cor}[teo]{Corollary}
\newtheorem{ackn}{Acknowledgments\!}

\theoremstyle{definition}

\theoremstyle{remark}

\newtheorem{rem}[teo]{Remark}

\numberwithin{equation}{section}

\def\SS{{{\mathbb S}}}

\def\RR{{\mathbb R}}

\def\RRR{{\mathrm R}}
\def\WWW{{\mathrm W}}

\def\TTT{{\mathrm{T}}}
\def\Ric{{\mathrm {Ric}}}

\def\SSS{{\mathrm S}}
\def\CCC{{\mathrm C}}

\def\dt{\frac{\partial\,}{\partial t}}

\title[Locally Conformally Flat Ancient Ricci Flows]{Locally Conformally Flat Ancient Ricci Flows}
\date{\today}

\author[Giovanni Catino]{Giovanni Catino}
\address[Giovanni Catino]{Dipartimento di Matematica, Politecnico di Milano, Piazza Leonardo da Vinci 32, Milano, Italy, 20133}
\email[G. Catino]{giovanni.catino@polimi.it}

\author[Carlo Mantegazza]{Carlo Mantegazza}
\address[Carlo Mantegazza]{Scuola Normale Superiore di Pisa, Piazza dei Cavalieri 7, Pisa, Italy, 56126}
\email[C. Mantegazza]{c.mantegazza@sns.it}

\author[Lorenzo Mazzieri]{Lorenzo Mazzieri}
\address[Lorenzo Mazzieri]{Scuola Normale Superiore di Pisa, Piazza dei Cavalieri 7, Pisa, Italy, 56126}
\email[L. Mazzieri]{l.mazzieri@sns.it}

\begin{document}

\begin{abstract} We show that any locally conformally flat ancient solution to the Ricci flow must be rotationally symmetric. As a by--product, we prove that any locally conformally flat Ricci soliton is a gradient soliton in the shrinking and steady cases as well as in the expanding case, provided the soliton has nonnegative curvature.
\end{abstract}

\maketitle

\begin{center}


\medskip

\centerline{\bf AMS subject classification:  53C24, 53C25, 53C44}

\end{center}


\section{Introduction}

In this paper, we study ancient solutions to the Ricci flow. We recall that a time--dependent metric $g(t)$ on a Riemannian manifold $M$ is a solution to the Ricci flow if it evolves by the equation
$$
\dt g(t) =-2 \,\Ric_{g(t)}\,.
$$
A solution is called {\em ancient} if it is defined for every negative time. Ancient solutions typically arise as the limit of a sequence of suitable blow--ups as the time approaches a singular time for the Ricci flow.  In dimension two there exists a compact, rotationally symmetric, ancient solution due to King~\cite{king}, Rosenau~\cite{rosenau} and Fateev, Onofri and Zamolodchikov~\cite{fatonozam}. In dimension three Perelman~\cite{perel2} constructed a compact, rotationally symmetric, ancient solution on the three sphere. In the non--rotationally symmetric case, the first construction is due to Fateev~\cite{fateev} in dimension three. Motivated by this construction, Bakas, Kong and Ni~\cite{bakkonni} produced high dimensional compact ancient solutions to the Ricci flow which are not rotationally symmetric.

In dimension two, Daskalopoulos, Hamilton and Sesum~\cite{dashamses2} have obtained a complete classification of all compact ancient solutions to the Ricci flow. Ni~\cite{leini5} showed that any compact ancient solution to the Ricci flow which is of type I, is $k$--noncollapsed, and has positive curvature operator has constant sectional curvature. In~\cite{brehuisin} Brendle, Huisken and Sinestrari proved that any compact ancient solution which satisfies a suitable pinching condition must have constant sectional curvature.

In this article, we show that any complete ancient solution to the Ricci flow in dimension $n\geq 4$ which is locally conformally flat along the flow must be rotationally symmetric.

\begin{teo}\label{mainancient} Let $(M^{n},g(t))$, $n\geq 4$, be a complete ancient solution to the Ricci flow which is locally conformally flat at every time. Then $(M^{n},g(t))$ is rotationally symmetric.
\end{teo}

The non--rotationally symmetric examples of Bakas, Kong and Ni show that the locally conformally flatness assumption cannot be removed. The proof of Theorem~\ref{mainancient} relies on a previous work of the first two authors~\cite{mancat1} about the behavior of the Weyl tensor under the Ricci flow, combined with a more recent result~\cite{catmanmaz} concerning the classification of Riemannian manifolds admitting a Codazzi tensor with exactly two distinct eigenvalues.

As a consequence of Theorem~\ref{mainancient}, we classify locally conformally flat Ricci solitons. We recall that a complete Riemannian manifold $(M^n ,g)$ is a {\em Ricci soliton} if there exists a vector field $X$ on $M^n$ such that 
$$
\Ric + \frac{1}{2}\mathcal{L}_{X}g \, = \, \lambda\, g
$$ 
for some constant $\lambda$. The Ricci soliton is called shrinking if $\lambda>0$, steady if $\lambda=0$ and  expanding if $\lambda<0$. If $X=\nabla f$ for some smooth function $f$, then the soliton is called a {\em gradient Ricci soliton}. It follows from the work of Perelman~\cite{perel1} (see~\cite{mantemin2}, for instance, for a direct proof) that any compact Ricci soliton is actually a gradient soliton. Moreover, Naber~\cite{na1} has showed that any shrinking Ricci soliton with bounded curvature has a gradient soliton structure. On the other hand, steady and expanding Ricci solitons which do not support a gradient structure were found by Lauret~\cite{lauret}, Baird and Danielo~\cite{bairdan}, Lott~\cite{lott} and Baird~\cite{baird}.

In this article we prove the following result.

\begin{teo}\label{mainsolitons} Let $(M^{n},g,X)$, $n\geq 4$, be a complete, locally conformally flat, shrinking or steady Ricci soliton. Then, it is a gradient Ricci soliton. The conclusion still holds in the expanding case, provided the soliton has nonnegative curvature operator.
\end{teo}

In particular, from the classification results of locally conformally flat gradient Ricci solitons in the shrinking case~\cite{caowang,mantemin2,nw2, pw, zhang}, in the steady case~\cite{caochen, mancat1}, as well as in the expanding case ~\cite{mancat1}, we obtain the following corollaries.

\begin{cor}\label{corshr} Let $(M^{n},g,X)$, $n\geq 4$, be a complete, locally conformally flat, shrinking Ricci soliton. Then, it is isometric to a quotient of $\SS^{n}$, $\RR\times \SS^{n-1}$ or $\RR^{n}$.
\end{cor}

\begin{cor}\label{corste} Let $(M^{n},g,X)$, $n\geq 4$, be a complete, locally conformally flat, steady Ricci soliton. Then, it is isometric to a quotient of $\RR^{n}$ or the Bryant soliton.
\end{cor}

\begin{cor}\label{corexp} Let $(M^{n},g,X)$, $n\geq 4$, be a complete, locally conformally flat, expanding Ricci soliton with nonnegative curvature operator. Then, it is a rotationally symmetric gradient expanding Ricci soliton.
\end{cor}

We notice that, rotationally symmetric gradient expanding Ricci solitons were constructed in~\cite{bry3, cao4, ilman6}.

\section{Notations and preliminaries}

The Riemann curvature
operator of a Riemannian manifold $(M^n,g)$ is defined 
as in~\cite{gahula} by
$$
\mathrm{Riem}(X,Y)Z=\nabla_{Y}\nabla_{X}Z-\nabla_{X}\nabla_{Y}Z+\nabla_{[X,Y]}Z\,.
$$ 
In a local coordinate system the components of the $(3,1)$--Riemann 
curvature tensor are given by
$\RRR^{l}_{ijk}\tfrac{\partial}{\partial
  x^{l}}=\mathrm{Riem}\big(\tfrac{\partial}{\partial
  x^{i}},\tfrac{\partial}{\partial
  x^{j}}\big)\tfrac{\partial}{\partial x^{k}}$ and we denote by
$\RRR_{ijkl}=g_{lm}\RRR^{m}_{ijk}$ its $(4,0)$--version.\\
With this choice, we have that the round sphere $\SS^n$ has positive
curvature, meaning that ${\mathrm{Riem}}(v,w,v,w)=\RRR_{ijkl}v^i w^j
v^k w^l >0$, for every couple $u$ and $v$ of non parallel vector
fields. The Ricci tensor is obtained by the contraction 
$\RRR_{ik}=g^{jl}\RRR_{ijkl}$ and $\RRR=g^{ik}\RRR_{ik}$ will 
denote the scalar curvature. The so called Weyl tensor is then 
defined by the following decomposition formula (see~\cite[Chapter~3,
Section~K]{gahula}) in dimension $n\geq 3$,
\begin{eqnarray}
\label{Weyl}
\WWW_{ijkl}  & = & \RRR_{ijkl} \, + \, \frac{\RRR}{(n-1)(n-2)} \,
(g_{ik}g_{jl}-g_{il}g_{jk}) \nonumber \\
& & - \, \frac{1}{n-2} \, (\RRR_{ik}g_{jl}-\RRR_{il}g_{jk}
+\RRR_{jl}g_{ik}-\RRR_{jk}g_{il}) \, .
\end{eqnarray}
The Weyl tensor shares the symmetries of the curvature
tensor. Moreover, as it can be easily seen by the formula above, all
of its contractions with the metric are zero.

In dimension three, $\WWW$ is identically zero on every Riemannian
manifold, whereas, when $n\geq 4$, the vanishing of the Weyl tensor is
a relevant condition, since it is  equivalent to the {\em locally
  conformally flatness} of $(M^n,g)$. We recall that this latter
condition means that around every point $p\in M^n$ there exists a
smooth function  $f$ defined in a open neighborhood $U_p$ of $p$, such
that the conformal deformation $\widetilde{g}$ of the original metric
$g$ defined by $\widetilde{g}_{ij}=e^fg_{ij}$ is flat. In particular,
the Riemann tensor associated to $\widetilde{g}$ is zero in $U_p$.

We also recall that in dimension $n=3$, {\em locally conformally
  flatness} is equivalent to the vanishing of the Cotton tensor
$$
\CCC_{ijk} \,\,= \,\, \nabla_{k} \RRR_{ij} - \nabla_{j} \RRR_{ik} \, - \, 
\frac{1}{2(n-1)} \, \big( \nabla_{k} \RRR \, g_{ij} - \nabla_{j} \RRR \,
g_{ik} \big)\,.
$$ 
By direct computation, we can see that the tensor $\CCC_{ijk}$
satisfies the following symmetries
\begin{equation}\label{CottonSym}
\CCC_{ijk}=-\CCC_{ikj},\,\quad\quad\CCC_{ijk}+\CCC_{jki}+\CCC_{kij}=0\,,
\end{equation}
moreover it is trace--free in any two indices, 
\begin{equation}\label{CottonTraces}
g^{ij}\CCC_{ijk}=g^{ik}\CCC_{ijk}=g^{jk}\CCC_{ijk}=0\,,
\end{equation}
by its skew--symmetry and Schur lemma. We notice that, for $n\geq 4$,
there holds
\begin{equation}\label{codweyl}
\,\nabla^l\WWW_{ijkl} \,\, = \,\, -\frac{n-3}{n-2} \, \CCC_{ijk} \, ,
\end{equation}
and we refer the reader to~\cite{besse} for the detailed
computation. It follows from this formula that, in every  dimension $n
\geq 3$, the vanishing of the Cotton tensor is a necessary condition
for a Riemannian manifold $(M^n, g)$ to be {\em locally conformally
  flat}. We also notice that the vanishing of the Cotton tensor can be
rephrased in terms of the so called Schouten tensor
$$
\SSS_{ij}=\RRR_{ij}-\frac{1}{2(n-1)}\RRR\, g_{ij} \,,
$$ 
by saying that $\SSS$ must satisfy the Codazzi equation
$$
(\nabla_{X}\SSS)\,Y \, = \, (\nabla_{Y}\SSS)\,X, \quad X,Y\in TM\,.
$$
Any symmetric two tensor satisfying this condition is called a Codazzi
tensor (see~\cite[Chapter~16]{besse} for a general overview on Codazzi
tensors). Hence, if $(M^n,g)$, $n\geq 3$, is a {\em locally
  conformally flat} manifold, then the Schouten tensor is a Codazzi
tensor.

\section{Proof of Theorem~\ref{mainancient}}

Let $(M^{n},g(t))$, $n\geq 4$, be a complete ancient solution to the Ricci flow. We assume that along the flow, the Weyl tensor remains identically zero. As it was observed in~\cite{mancat1}, this condition implies a strong rigidity on the eigenvalues of the Ricci tensor. More precisely, one has the following result.

\begin{lemma}[{\cite[Corollary 1.2]{mancat1}}]\label{lem_eigen}
Let $(M^n,g)$, $n\geq 4$, be a solution to the Ricci flow such that the Weyl tensor remains identically zero at every time. Then, at every point, either the Ricci tensor is proportional to the metric or it has an eigenvalue of multiplicity $(n-1)$ and another of multiplicity $1$.
\end{lemma}

By the results in~\cite{chen2} and~\cite{zhang}, which generalize the
well--known Hamilton--Ivey curvature estimate, we know that every
complete ancient solution $g(t)$ to the Ricci flow whose Weyl tensor
is identically zero for all times, is forced to have nonnegative
curvature operator for every time 
$t$. Moreover, by Hamilton's strong maximum principle for systems
in~\cite{hamilton2}, we have that either
the metric has strictly positive curvature operator or it splits a
line. By Theorem~1.167 in~\cite{besse} a Riemannian product $(\RR\times
N^{n-1}, ds\times h)$  is locally
conformally flat if and only if the manifold $(N^{n-1},h)$ has constant
curvature, hence, either one of the following possibilities holds:
$(M^{n}, g(t))$ is flat or it 
is a quotient of a rescaling of $\RR\times\SS^{n-1}$ or it has positive curvature
operator. Since the first two cases satisfy the conclusion of the
theorem, from now on we assume that $(M^{n}, g(t))$ is a complete,
locally conformally flat, ancient solution to the Ricci flow with
positive curvature operator.

As we have seen in the previous section, the relation~\eqref{codweyl} implies that the Cotton tensor is identically
zero, hence the Schouten tensor
$$
\SSS_{ij}=\RRR_{ij}-\frac{1}{2(n-1)}\RRR\, g_{ij}
$$ 
is a Codazzi tensor. Moreover, from Lemma~\ref{lem_eigen}, we know that, at every point, either the metric is Einstein or the Ricci tensor (and so the Schouten tensor) has two distinct eigenvalues of multiplicity $1$ and $(n-1)$, respectively.\\
Now, it was proved by Bando~\cite{bando} that solutions to the Ricci flow are real analytic. To be precise, Bando showed that any Ricci flow solutions is real analytic if $M^{n}$ is compact or if it is complete with uniformly bounded curvature. This result was recently improved by Kotschwar~\cite{kotschwar2}, who showed a local version of Bando's result. It follows that if the metric is Einstein in some open subset of $M^n$, then it is Einstein everywhere, and by conformally flatness and positivity of the curvature the manifold $(M^n,g)$ must be isometric to a quotient of $\SS^{n}$. Thus, either $(M^{n},g(t))$ has constant positive sectional curvature or the Schouten tensor has an eigenvalue of multiplicity $1$ and a different one of multiplicity $(n-1)$ at every point of some open dense subset $U$ of $M^{n}$. In the latter case, we apply to such open set $U$ the following classification result of Riemannian manifolds admitting a Codazzi tensor with two distinct eigenvalues.

\begin{lemma}[{\cite[Theorem 2.1]{catmanmaz}}]\label{lem_codazzi}
Let $\TTT$ be a Codazzi tensor on the Riemannian manifold $(U,g)$, with $n\geq 3$. Suppose that at every point of $U$, the tensor $\TTT$ has exactly two distinct eigenvalues $\rho$ and $\sigma$ of multiplicity $1$ and $n-1$, respectively. Finally, we let $W=\{ p\in U \, \big| \, d\sigma(p)\neq 0 \}$. Then, we have that 

\begin{enumerate}

\item The closed set $\overline{W}=W\cup \partial W$ with the metric $g|_{\overline{W}}$ is locally isometric to the warped product of some $(n-1)$--dimensional Riemannian manifold on an interval of $\RR$ and $\sigma$ is constant along the "leaves" of the warped product.

\item The boundary of $W$, if present, is given by the disjoint union of connected totally geodesic hypersurfaces where $\sigma$ is constant. 

\item Each connected component of the complement of $\overline{W}$ in $U$, if present, has $\sigma$ constant and it is foliated by totally geodesic hypersurfaces.
\end{enumerate}
\end{lemma}

We are going to to show that, under our assumptions, case $(3)$ cannot occurs and $W = U$. In fact, if a connected component of the complement of $\overline{W}$ in $U$ is present, actually, the proof of this lemma 
(see also~\cite[Proposition 16.11]{besse}) shows that the totally geodesic foliation of $(U,g(t))$ is {\em integrable}. Then, since the manifold $(U,g(t))$ 
has nonnegative sectional curvature, it follows from~\cite[Corollary 2]{abe}, that such component must split a flat factor. 
This is clearly in contradiction with the positivity of the curvature, hence it must be $U=W$ in this lemma and 
$(U,g(t))$ is locally a warped product of some $(n-1)$--dimensional manifold on an interval of $\RR$.\\
Since $(U,g(t))$ is locally conformally flat with positive curvature operator, we have that the $(n-1)$--dimensional fibers of the warped product are 
isometric to $\SS^{n-1}$ and the metric is rotationally symmetric. 
By the density of $U$ in $M^n$, this conclusion clearly holds for the whole $(M^n,g)$. This concludes the proof of Theorem~\ref{mainancient}.

\ 
 
\begin{rem}\label{rem1}
We would like to notice that, the same argument shows that the conclusion of Theorem~\ref{mainancient} still holds if one consider a Ricci flow solution $(M^{n},g(t))$, $n\geq 4$, defined on some time interval $I\subseteq\RR$, which is locally conformally flat with nonnegative curvature operator for every $t\in I$.
\end{rem}

\section{Proof of Theorem~\ref{mainsolitons}}

Now we turn our attention to the classification of locally conformally flat Ricci solitons. Let $(M^{n},g, X)$ be a complete, locally conformally flat shrinking or steady Ricci soliton. In particular, it generates a self--similar ancient solution $g(t)$ to Ricci flow (see~\cite{zhang2}) which is locally conformally flat at every time $t$. Hence, Theorem~\ref{mainancient} implies that the manifold is rotationally symmetric with nonnegative curvature operator. As we observed in Remark~\ref{rem1}, the conclusion still holds if we consider an expanding Ricci soliton with nonnegative curvature operator. To prove Theorem~\ref{mainsolitons} we then apply the following result.

\begin{lemma} [{\cite[Proposition 2.6]{mancat1}}] 
Let $(M^n,g, X)$ be a complete, locally warped, locally conformally flat Ricci soliton with nonnegative Ricci tensor, then it is a gradient Ricci soliton with a potential function $f:M^n\to\RR$ (hence, $X=\nabla f$) depending only on the $r$ variable of the warping interval.
\end{lemma}

\medskip

\begin{ackn} The authors are partially supported by the Italian projects FIRB--IDEAS ``Analysis and Beyond'' and GNAMPA projects ``Flussi geometrici e soluzioni autosimilari'' and ``Equazioni di evoluzione geometriche e strutture di tipo Einstein''.
\end{ackn}

\bibliographystyle{amsplain}
\bibliography{biblio}

\end{document}